\documentclass[11pt]{amsart}

\usepackage{amsfonts}
\usepackage{amsmath}
\usepackage{amssymb}
\usepackage{amsthm}
\usepackage{enumitem}
\usepackage{etoolbox}

\usepackage{array}   
\newcolumntype{L}{>{$}c<{$}} 

\usepackage{tikz}
\usetikzlibrary{cd}
\let\amsamp=&

\usepackage{graphicx}
\usepackage{hyperref}
\usepackage{pdflscape}

\usepackage{algorithm2e}

\usepackage{quiver}

\usepackage[capitalize]{cleveref}

\newcommand{\tensor}{\otimes}

\DeclareMathOperator{\Hom}{Hom}
\DeclareMathOperator{\Sq}{Sq}

\renewcommand{\epsilon}{\varepsilon}
\renewcommand{\phi}{\varphi}

\newtheorem{conj}{Conjecture}[section]
\newtheorem{lemma}[conj]{Lemma}
\newtheorem{prop}[conj]{Proposition}
\newtheorem{thm}[conj]{Theorem}

\theoremstyle{definition}
\newtheorem*{convention}{Conventions}
\newtheorem{defn}[conj]{Definition}
\newtheorem{example}[conj]{Example}
\newtheorem{rmk}[conj]{Remark}

\Crefname{rmk}{Remark}{Remarks}
\Crefname{rmkenumi}{Remark}{Remarks}
\AtBeginEnvironment{rmk}{%
    \crefalias{enumi}{rmkenumi}%
    \setlist[enumerate,1]{
        label={\textit{(\roman*)}},
        ref={\thermk(\roman*)}
    }%
}

\begin{document}

\author{Joey Beauvais-Feisthauer}

\title[Automated Adams differential computation]{Automated differential computation in the Adams Spectral Sequence}

\address{Department of Mathematics, Wayne State University, Detroit, MI 48202}

\email{joeybf@wayne.edu}

\thanks{The author was partially supported by National Science Foundation grant DMS-1904241.}



\begin{abstract}
    We describe an algorithm for the automated deduction of many $d_2$ differentials in the Adams spectral sequence. We discuss our implementation and the results of the computation.
\end{abstract}

\maketitle

\section{Introduction}

A central problem in homotopy theory is the computation of stable homotopy groups. One of our most efficient tools is the Adams spectral sequence, which uses information about certain Ext groups to compute stable homotopy groups. In our case of interest, its $E_2$ page is a differential graded algebra isomorphic to the cohomology of the Steenrod algebra, and it converges to the stable homotopy groups of spheres.

The determination of the additive and multiplicative structure of this differential graded algebra is a hard problem, but manageable in practice. Being defined entirely in terms of homological algebra, its calculation can be automated using computers. For instance, Bruner and Rognes~\cite{bruner2021data} have computed it up to total degree 184, which is the data that we use for our results in this paper. They have since then extended their computation to total degree 200~\cite{bruner2022data}, and the author is currently participating in a similar computation up to stem 256~\cite{bf2022}.

While we have a good grasp of the algebraic structure of its $E_2$ page, the main challenge when working with the Adams spectral sequence in general is the computation of its differentials. In particular, the first obstacle is computing the differentials on the $E_2$ page. A wide variety of ad hoc techniques have been used to compute some of their values~\cite{isaksen2020more}. We have developed a technique to systematize this procedure.

We introduce the new results that our algorithm found in Section~\ref{sec:results}. We then give a brief exposition of the theory behind our algorithm in Sections~\ref{sec:affine}~and~\ref{sec:leibniz}, followed by a breakdown of the algorithm in Section~\ref{sec:algo}.

\subsection*{Related work}

We should mention that some of these results were already present in unpublished work of Bruner from around 2005~\cite{bru05}. However, we believe that our method is more conceptual and generalizes better to other spectral sequences. Furthermore, we apply this algorithm to a larger range, which necessarily gives us stronger results.

Recently, Dexter Chua has developed a tool to compute $d_2$ differentials algorithmically, directly from homological algebra~\cite{chua2021adams}. In practice we have found that, using Chua's algorithm, computing the differentials on the $E_2$ page is no harder than computing its product structure. This arguably means that, similarly to the computation of the additive and multiplicative structure of the $E_2$ page, the computation of its differentials can be considered ``solved''.

However, we believe our method of propagation remains relevant. Firstly, the algorithm allows us to not only determine many $d_2$ differentials but explore their interrelationship. Noticing this qualitative behavior is what led us to formulate Conjecture~\ref{conj:big}. Secondly, although we only discuss the results that our procedure gave while examining the Adams spectral sequence, and more specifically its $E_2$ page, this same procedure is readily applicable to any other differential graded algebra. It has already been used to study the $d_3$ differentials of the Adams spectral sequence, which have been found to behave very differently from the $d_2$ differentials. The code can also be easily adapted to work with any spectral sequence equipped with a product structure, such as those of May type.

\begin{convention}\mbox{}
    \begin{enumerate}
        \item We index the Adams spectral sequence using bidegrees $(n,s)$ where $n$ is the stem and $s$ is the homological degree, i.e. the usual Cartesian coordinates on an Adams chart. We use capital latin letters $A, B, \ldots$ to refer to individual bidegrees, e.g. $A = (20,4)$. Addition of bidegrees is componentwise addition. Also, $A'$ denotes the bidegree of the output of the differential on a given bidegree $A$. In other words, if $A = (n,s)$, then $A' = (n-1, s+2)$. As an example, if we are in the $E_2$ page of the Adams spectral sequence, and $A = (20,4)$, then $A' = (20-1,4+2) = (19,6)$.
        
        Notice that, for any bidegrees $A$ and $B$,
        $$ A' + B = A + B' = (A + B)'. $$
    
        \item We denote the differential bigraded algebra by $E$, its homogeneous component in bidegree $A$ by $E^A$, and the differential on $E^A$ (as a linear map) by $d^A$. The symbol $D^A$ always denotes an \emph{affine} subspace of $\Hom(E^A, E^{A'})$ (see Section~\ref{sec:affine}), and we assume that $d^A$ belongs to $D^A$ unless otherwise noted.
    \end{enumerate}
    
\end{convention}

\section{Results}\label{sec:results}
Here are some results that our algorithm has yielded. The data we used for the computations was generated in 2020 by Bruner and Rognes, using Bruner's software~\cite{bruner2018ext}. It is a complete description of the cohomology of the Steenrod algebra through total degree 184, along with its primary multiplicative structure. As mentioned in the introduction, Bruner and Rognes have since extended their computation through total degree 200. We have preliminary results on the application of our algorithm to the larger computation of Bruner and Rognes. For expository precision, we restrict our discussion only to the smaller computation of Bruner and Rognes.

It is well-known that $d_2(h_4) = h_0 h_3^2$. Adams gives the historically first proof~\cite{MR141119}, but this differential can also be computed by other arguments. For example, it can be derived by recognizing $h_4 = \Sq^0(h_3)$ and using Bruner's theorem on the interaction between differentials and algebraic squaring operations~\cite{MR836132}. The more concretely-minded reader might prefer to read this differential off a machine computation, such as the one in~\cite{chua2021adams}.

Similarly, it is known that $d_2(\Delta^2 d_0^2) = d_0 j m$. This differential can again be obtained by several arguments. For example, it follows by comparison with the Adams spectral sequence for $\mathrm{tmf}$~\cite{MR4284897}, but it can also be read off of a machine computation.

Given only the knowledge of $d_2(h_4) = h_0 h_3^2$, our algorithm is able to infer the value of the differentials on a large portion of the $E_2$ page. Given only the values of $d_2(h_4)$ and $d_2(\Delta^2 d_0^2)$, we were able to compute over 95\% of the differentials up to stem 140. In particular, we found the first known proofs for the values of the differentials on the following elements:
\[\begin{tabular}{LLLL}
\hline
\text{stem} & \text{filtration} & x & d_2(x) \\
\hline
104 & 18 & \Delta^4 h_1 h_3 & M P \Delta^2 h_1^2 + M \Delta h_2^2 d_0^2 \\
107 & 13 & \Delta D_1 n & 0 \\
107 & 18 & M \Delta^2 d_0 & M \Delta h_2^2 d_0 e_0 \\
108 & 14 & \Delta^2 \Delta_1 h_1 h_3 & M P \Delta_1 h_1^2 \\
108 & 15 & M \Delta^2 h_4 & 0 \\
109 & 12 & x_{109,12} & 0 \\
109 & 14 & \Delta^2 A' & \Delta g^2 g_2 \\
109 & 14 & \Delta^2 (A+A') & M \Delta^2 h_0 h_4 \\
110 & 18 & \Delta^4 h_3^2 & h_0^3 x_{109,17} \\
110 & 18 & M \Delta^2 e_0 & M \Delta h_2^2 e_0^2 \\
\hline
\end{tabular}\]
See \cite{isaksen2020more} for an explanation of the notation.

\begin{thm}\mbox{}\label{thm:main}
    \begin{enumerate}
        \item The product structure of the Adams $E_2$ page implies that all elements represented by a black dot in Figure~\ref{fig:diff_baseline} do not support a $d_2$ differential.
        
        \item The product structure of the Adams $E_2$ page, together with the fact that $d_2(h_4) = h_0 h_3^2$, implies the existence of all cyan differentials in Figure~\ref{fig:diff_h4_7616}.
        
        \item The product structure of the Adams $E_2$ page, together with the facts that $d_2(h_4) = h_0 h_3^2$ and $d_2(\Delta^2 d_0^2) = d_0 j m$, implies:
        \begin{enumerate}
            \item that all elements represented by a black dot in Figure~\ref{fig:diff_h4_7616} do not support a $d_2$ differential.
            
            \item the existence of all cyan differentials and all magenta differentials.
        \end{enumerate}
    \end{enumerate}
\end{thm}

\begin{landscape}\centering
\vspace*{\fill}
\begin{figure}[htpb]
    \centering
    \includegraphics[
    width=1.6\textwidth]{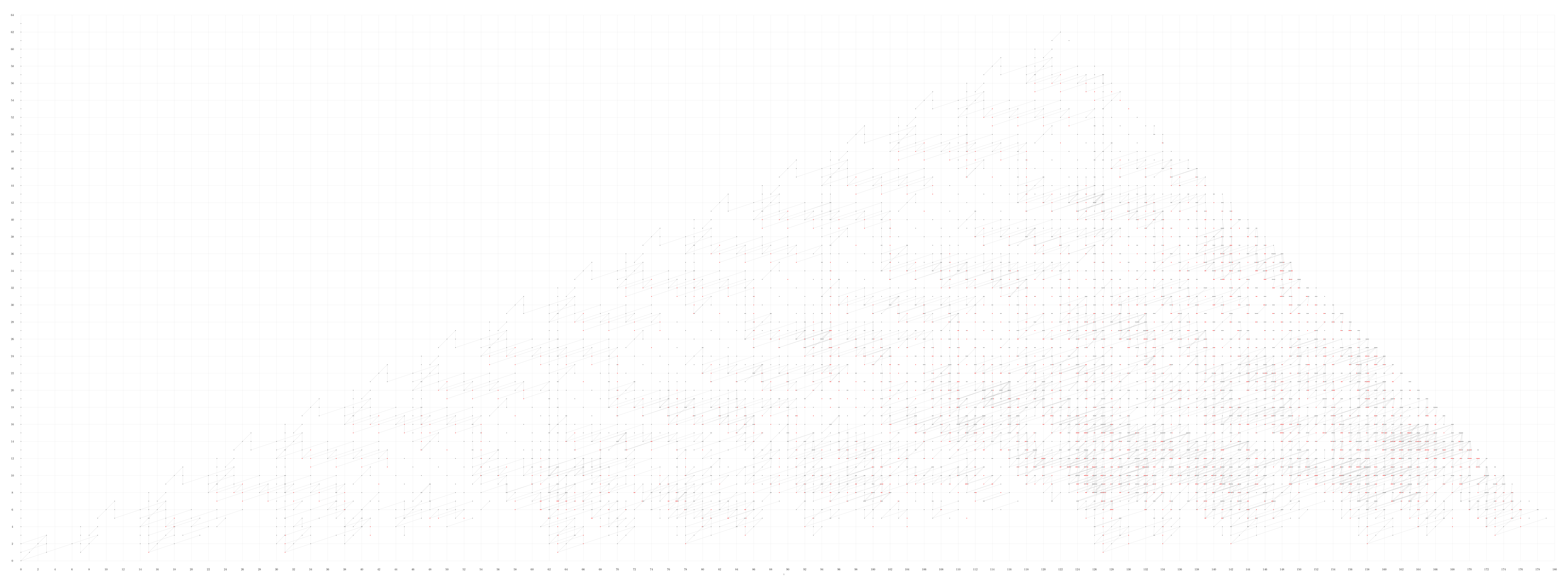}
    \caption{Elements not supporting a differential (See Theorem~\ref{thm:main})}
    \label{fig:diff_baseline}
\end{figure}
\vfill
\end{landscape}

\begin{landscape}\centering
\vspace*{\fill}
\begin{figure}[htpb]
    \centering
    \includegraphics[
    width=1.6\textwidth]{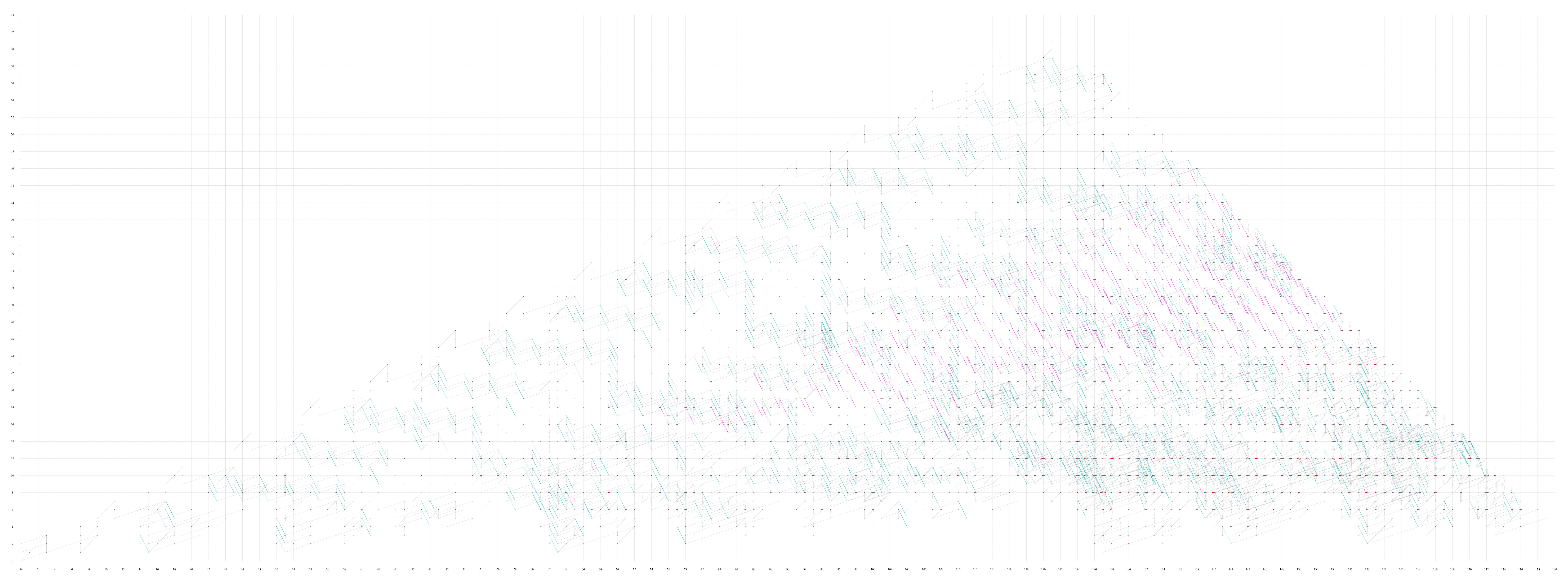}
    \caption{Differentials derived by the algorithm (See Theorem~\ref{thm:main})}
    \label{fig:diff_h4_7616}
\end{figure}
\vfill
\end{landscape}

Although our tool was unable to derive all differentials, it highlighted which ones could not be calculated by routine arguments. This led to a careful manual double-check of these critical differentials, which uncovered a mistake in the computation of the $d_2$ differentials in the 85-stem~\cite{isaksen2020more}.

Here is an example deduction that our tool is able to perform. The proposition relates two statements that are known to be true. However, its proof might be of value, both because we believe it cannot be found in the literature, and because it was found in an entirely automated fashion. We hope that this proof foreshadows more extensive machine deduction in the future of stable homotopy theory. We emphasize that our written proof merely transcribes in human language some of the output of our algorithm.

\begin{prop}
The element $h_4$ supports a differential if and only if $e_0$ supports a differential.
\end{prop}
\begin{proof}
Using the relation $h_4 \cdot i = 0$, we get
$$ d_2(h_4) i + h_4 d_2(i) = 0. $$
The bidegree containing $d_2(h_4)$ does not contain any $i$-torsion elements, and similarly the bidegree containing $d_2(i)$ does not contain any $h_4$-torsion elements. Therefore, $h_4$ supports a differential if and only if $i$ also supports a differential.

Using the relation $f_0 \cdot i = h_1^2 \cdot \Delta h_1 d_0$, we get
$$ d_2(f_0) i + f_0 d_2(i) = d_2(h_1^2 \cdot \Delta h_1 d_0). $$
We can see that $h_1^2 \cdot \Delta h_1 d_0$ is divisible by $h_1^2$, and the only nonzero element in the bidegree containing its differential is $h_1 \cdot d_0 \cdot P e_0$, which is not divisible by $h_1^2$. Therefore, the right-hand side of the equation must be zero. Multiplication by $i$ on the bidegree containing $d_2(f_0)$ is an isomorphism, and the same is true of multiplication by $f_0$ on the bidegree containing $d_2(i)$. Therefore, $i$ supports a differential if and only if $f_0$ also supports a differential.

Finally, using the relation $h_1 e_0 = h_0 f_0$, a similar argument gives us that $f_0$ supports a differential if and only if $e_0$ also does.
\end{proof}

Our procedure is able to determine most of the $d_2$ differentials in relatively low stems, with a few notable exceptions.

For example, inspection of Figure~\ref{fig:diff_h4_7616} shows that the differential of the element $D_1$ in bidegree $(52,5)$ is not computed by our method. This is due to the fact that the only possible nonzero value of $d_2(D_1)$ is $h_2^2 h_5 d_0$. The element $h_2^2 h_5 d_0$ annihilates the entire $E_2$ in our range, so the multiplicative structure alone cannot distinguish between it and 0. However, more recent calculations have shown that in fact $h_2^2 h_5 d_0$ has a nonzero $h_7$ multiple in bidegree $(178, 8)$. Unfortunately, this is still not enough for our algorithm to compute the differential. More generally, a result of Don Davis~\cite{davis1981infinite} shows that any element of the $E_2$ page has nonzero $h_i$ multiples for all large enough $i$. This result guarantees that, in principle, any given element will have arbitrarily many nonzero products in a large enough range. Heuristically, as the number of nonzero products increases, there is also an increasing number of relations between them, which would allow our algorithm to backtrack from these later differentials to infer results in relatively low stems.

This heuristic argument leads us to formulate the following conjecture.
\begin{conj}\label{conj:big}
    The space of derivations on the Adams $E_2$ page for the sphere is 1-dimensional. In other words, there is only one nontrivial set of differentials on $E_2$ which is compatible with its algebraic structure.
\end{conj}

\subsection*{Further research}
Our tool does not currently enforce that $d^2 = 0$. It is known that enforcing this requirement would restrict the possible differentials even further. The first example happens in bidegree $(125,10)$: our tool cannot rule out a certain matrix for the differential on that bidegree, even though there is no allowable differential on $(124,12)$ which would yield 0 when composed with this matrix.

We chose not to add this condition because it would require us to work with non-affine subspaces. It is very efficient to encode an affine subspace as an offset vector and a list of basis elements, and every computation in the algorithm can be reduced to operations on those basis elements. Keeping track of a non-affine subspace would increase the complexity of the code, the memory usage and the execution time by a large factor.

In~\cite{bru05}, Bruner used algebraic squaring operations to compute differentials. We have not attempted to implement this functionality, even though we strongly suspect that using these operations would give us a number of new differentials.

\section{Affine subpaces}\label{sec:affine}
Let $V$, $W$, and $U$ be vector spaces and $f: V \to W$ be a linear map. Recall that the map $f$ gives rise to the linear maps
\begin{align*}
    f^* &: \Hom(W, U) \to \Hom(V, U), \\
    f_* &: \Hom(U, V) \to \Hom(U, W),
\end{align*}
given by pre- and post-composing with $f$ respectively.
These operations are themselves linear maps
\begin{align*}
    (-)^*: \Hom(V, W) &\to \Hom(\Hom(W, U), \Hom(V, U)), \\
    (-)_*: \Hom(V, W) &\to \Hom(\Hom(U, V), \Hom(U, W)).
\end{align*}
In categorical language, these maps are the components of the internal Yoneda and co-Yoneda embeddings respectively, where we are taking the $U$-component in the codomain.

\begin{rmk} \label{upper_star}
It makes sense to apply $(-)^*$ and $(-)_*$ to any set $D \subseteq \Hom(V, W)$. We will denote their images by $D^*$ and $D_*$ respectively.
\end{rmk}

Our algorithm manipulates affine subspaces extensively. Since they are less commonly studied than linear subspaces, we review some theory in this subsection.

\begin{defn}
    Let $V$ be a vector space. An \emph{affine subspace} of $V$ is a subset of the form $v + U$ where $U$ is a linear subspace of $V$, and $v$ is a vector. We call $v$ the \emph{offset vector} or simply \emph{offset}, and $U$ the \emph{linear part}. The dimension of the affine subspace is the dimension of the linear part $U$. We say an affine subspace is \emph{linear} if it contains the origin.
\end{defn}

\begin{rmk} \label{affine_remark}\mbox{}
    \begin{enumerate}
        \item Affine subspaces don't have a canonical offset in general. If $v + U$ is an affine subspace of $V$, we have $v + U = (v + u) + U$ for any element $u$ in $U$. In fact, $v + U = w + U$ if and only if $v - w$ belongs to $U$.
        
        \item As a special case, vectors themselves are affine subspaces, with linear part equal to the 0 subspace. We often identify without mention an affine subspace with its offset vector in the 0-dimensional case.
        
        \item \label{rem:affine_not_empty} Affine spaces are never empty, since they contain at least the offset vector.
    \end{enumerate}
\end{rmk}

Because of \cref{rem:affine_not_empty}, affine spaces are not closed under arbitrary intersections. However, the next best thing is true.
\begin{lemma}\label{lem:affine_inter}
    Let $V$ be a vector space, and let $\{D_i\}_{i \in I}$ be an arbitrary family of affine subspaces. If their intersection $D = \bigcap_{i \in I} D_i$ is nonempty, then $D$ is an affine subspace.
\end{lemma}
\begin{proof}
    Let $v$ belong to the intersection. We can take $v$ as the offset vector for all subspaces, so
    $$ D_i = v + U_i $$
    for some family $U_i$ of \emph{linear} subspaces of $V$. The intersection is then exactly the affine subspace
    \[ D = v + \bigcap_{i \in I} U_i. \qedhere \]
\end{proof}

\begin{defn}
    Let $V$ be a vector space, and $B \subseteq V$ be a nonempty subset. The \emph{affine span} of $B$ is the smallest affine subspace of $V$ containing $B$.
\end{defn}

To see that the affine span exists, it suffices to take the intersection of all affine spaces containing $B$, which is affine by Lemma~\ref{lem:affine_inter}.
Notice that if $B$ contains the origin, the affine span coincides with the usual linear span.

\begin{lemma} \label{lem:concrete_span}
    Let $V$ be a vector space, let $B \subseteq V$ be a nonempty subset and let $b$ be an arbitrary element of $B$. The affine span of $B$ is exactly the set
    $$ b + \langle\{b' - b \mid b' \in B\}\rangle. $$
\end{lemma}
\begin{proof}
    This set is affine and contains all elements of $B$. Conversely, let $v + U$ be an affine subspace containing $B$. Since $b$ belongs to this space, we may choose it as the offset vector. Let $b'$ be any element of $B$. Since $b'$ also belongs to this space, we have
    $$ b' + U = b + U, $$
    which implies that $b' - b$ belongs to $U$.
\end{proof}

\begin{defn}
    Let $V$ and $W$ be vector spaces, and let $D = v + U$ and $D' = v' + U'$ be affine subspaces of $V$ and $W$ respectively.
    \begin{enumerate}
        \item The \emph{direct sum} $D \oplus D'$ is the set
        $$ \{ (v, w) \mid v \in D, w \in D' \}. $$
        \item The \emph{tensor product} $D \tensor D'$ is the affine span of the set
        $$ \{ v \tensor w \mid v \in D, w \in D' \}. $$
    \end{enumerate}
\end{defn}

In particular, if $D = v$ is 0-dimensional and $D' = U'$ is linear, then
$$ D \tensor D' = \{ v \tensor w \mid w \in U'\} = v \tensor U'. $$
Notice that the tensor product is a linear subspace in this special case.

Affine subspaces, like linear subspaces, are closed under some common operations.

\begin{lemma}\label{lem:basic_affine}
    Let $V$ and $W$ be vector spaces. Let $D = v + U$ be an affine subspace of $V$ and $D' = v' + U'$ an affine subspace of $W$.
    \begin{enumerate}
        \item We have the equality of sets
        $$ D \oplus D' = (v, v') + (U \oplus U'). $$
        
        \item We have the equality of sets
        $$ D \tensor D' = (v \tensor v') + (v \tensor U' + U \tensor v' + U \tensor U'). $$
        In particular, if $D' = v'$ is 0-dimensional, then we have
        $$ D \tensor D' = (v \tensor v') + (U \tensor v'). $$
        
        \item Let $f: V \to W$ be a linear map. We have the equality of sets
        $$ f(D) = f(v) + f(U). $$

        \item Let $f: V \to W$ be a linear map. Assume that $D'$ intersects $f(V)$ non-trivially, and let $x$ be any vector such that $f(x)$ is in $D'$. We have the equality of sets
        $$ f^{-1}(D') = x + f^{-1}(U'). $$
    \end{enumerate}
\end{lemma}
Notice that in Lemma 3.8(4) we had to impose an extra condition, namely that $D'$ and $f(V)$ intersect non-trivially. This is necessary, because otherwise $f^{-1}(D')$ would be empty (see \cref{rem:affine_not_empty}).
\begin{proof}\mbox{}
    \begin{enumerate}
        \item \mbox{}
        \vspace{-\baselineskip}
        \vspace{-\abovedisplayskip}
        \begin{align*}
            D \oplus D' &= \{(x, x') \mid x \in D, x' \in D'\} \\
            &= \{(v + u, v' + u') \mid u \in U, u' \in U'\} \\
            &= \{(v, v') + (u, u') \mid u \in U, u' \in U'\} \\
            &= (v, v') + (U \oplus U').
        \end{align*} \mbox{}
    
        \item Let $B = \{a \tensor b \mid a \in D, b \in D'\} = \{(v + u) \tensor (v' + u') \mid u \in U, u' \in U'\}$. Note that $v \tensor v'$ belongs to $B$. Let
        $$ K = \langle \{(v + u) \tensor (v' + u') - v \tensor v' \mid u \in U, u' \in U'\} \rangle. $$
        By Lemma~\ref{lem:concrete_span}, it suffices to show that
        \begin{align*}
            K = v \tensor U' + U \tensor v' + U \tensor U'.
        \end{align*}
        Rewriting $K$ as
        $$ K = \langle \{v \tensor u' + u \tensor v' + u \tensor u' \mid u \in U, u' \in U'\} \rangle, $$
        we have the inclusion
        $$ K \subseteq v \tensor U' + U \tensor v' + U \tensor U'. $$
        Setting $u = 0$ shows that $K$ contains all elements of the form $v \tensor u'$, and similarly, setting $u' = 0$ shows that $K$ contains all elements of the form $u \tensor v'$. Moreover, for any $u$ in $U$ and $u'$ in $U'$,
        $$ u \tensor u' = (v \tensor u' + u \tensor v' + u \tensor u') - (v \tensor u') - (u \tensor v') \in K. $$
        This proves the reverse inclusion
        $$ K \supseteq v \tensor U' + U \tensor v' + U \tensor U'. $$ \mbox{}
        
        \item\mbox{}
        \vspace{-\baselineskip}
        \vspace{-\abovedisplayskip}
        \begin{align*}
            f(D) &= \{ f(x) \mid x \in D \} \\
            &= \{ f(v + u) \mid v + u \in v + U \} \\
            &= \{ f(v) + f(u) \mid v + u \in v + U \} \\
            &= f(v) + \{ f(u) \mid u \in U \} \\
            &= f(v) + f(U).
        \end{align*} \mbox{}
        
        \item Since $f(x)$ belongs to $D'$, we have $D' = f(x) + U'$. Therefore,
        \begin{align*}
            f^{-1}(D') &= \{ u \in V \mid f(u) \in D' \} \\
            &= \{ x + u' \in V \mid f(x + u') \in f(x) + U' \} \\
            &= x + \{ u' \in V \mid f(u') \in U' \} \\
            &= x + f^{-1}(U'). \qedhere
        \end{align*}
    \end{enumerate}
\end{proof}

\section{Leibniz rule}\label{sec:leibniz}

As a differential algebra, the product on the $E_2$ page satisfies the Leibniz rule. Consider two elements $x$ and $y$ belonging to the $E_2$ page. The Leibniz rule
$$ d_2(xy) = d_2(x)y + xd_2(y) $$
gives a linear relation between $d_2(x)$, $d_2(y)$ and $d_2(xy)$, such that information on any of them gives information on the possible values of the others.

We wish to restate the Leibniz rule in an element-free way. This formulation lets us manipulate subspaces directly instead of individual elements, which is crucial for our algorithm. Let $\mu$ be the multiplication
$$ \mu: E \tensor E \to E. $$
Suppose that $x$ is in bidegree $A$ and $y$ is in bidegree $B$. The Leibniz rule can be expressed in diagram form as follows:

\begin{equation}\label{diag:leibniz_normal}
    \begin{tikzcd}
    	{E^A \tensor E^B} &&& {E^{A+B}} \\
    	\\
    	{(E^{A'} \tensor E^B) \oplus (E^A \tensor E^{B'})} &&& {E^{(A+B)'}}
    	\arrow["{d^{A+B}}", from=1-4, to=3-4]
    	\arrow["\mu", from=1-1, to=1-4]
    	\arrow["{\begin{bmatrix} \mu \\ \mu \end{bmatrix}}", from=3-1, to=3-4]
    	\arrow["{\begin{bmatrix} d^A \tensor 1 \amsamp 1 \tensor d^B \end{bmatrix}}", from=1-1, to=3-1]
    \end{tikzcd}
\end{equation}
We annotate the diagram with the appropriate bidegrees to make it clear that the various $d$'s that appear should be considered as independent variables. The core of the algorithm will be to, in essence, fix two of the differentials and use the commutation of the diagrams to infer information about the third.

This diagram implies in particular the following proposition.
\begin{prop} \label{prop:leibniz_easy}
Let $D^A$ and $D^B$ be affine subspaces of $\Hom(E^A, E^{A'})$ and $\Hom(E^B, E^{B'})$ respectively. If the differential $d^A$ belongs to $D^A$ and the differential $d^B$ belongs to $D^B$, then the differential $d^{A+B}$ belongs to the affine subspace
$$ S(A,B) = (\mu^*)^{-1}(\mu_*(D^A \tensor 1) + \mu_*(1 \tensor D^B)) $$
of $\Hom(E^{A+B}, E^{(A+B)'})$. Here we have
\begin{align*}
    D^A \tensor 1 &\subseteq \Hom(E^A \tensor E^B, E^{A'} \tensor E^B), \\
    \mu_*(D^A \tensor 1) &\subseteq \Hom(E^A \tensor E^B, E^{(A+B)'}), \\
    1 \tensor D^B &\subseteq \Hom(E^A \tensor E^B, E^A \tensor E^{B'}), \\
    \mu_*(1 \tensor D^B) &\subseteq \Hom(E^A \tensor E^B, E^{(A+B)'}).
\end{align*}
\end{prop}
\begin{proof}
Recall that vector addition is a linear map
$$ +: V \oplus V \to V. $$
With this notation, $S(A,B)$ can be written
$$ S(A,B) = (\mu^*)^{-1}(+(\mu_*(D^A \tensor 1) \oplus \mu_*(1 \tensor D^B))). $$
Repeated application of Lemma~\ref{lem:basic_affine} to this description shows that $S(A,B)$ is indeed an affine subspace of $\Hom(E^{A+B}, E^{(A+B)'})$.

To see that it contains $d^{A+B}$, we note that equating the two composites in \eqref{diag:leibniz_normal} gives
$$ (d^A \tensor 1) \circ \mu + (1 \tensor d^B) \circ \mu = \mu \circ d^{A+B}. $$
By assumption, this map is contained in the space
$$ \mu_*(D^A \tensor 1) + \mu_*(1 \tensor D^B). $$
Taking preimages under $\mu^*$ gives the desired result.
\end{proof}

Notice that there is no mention of any element in Proposition~\ref{prop:leibniz_easy}, only linear maps and affine subspaces. Therefore, we have achieved our goal of restating the Leibniz rule in an element-free way. Next, we aim to modify this diagram so that one of the maps is $d^A$ and the others are expressed in terms of $d^B, d^{A+B}$ and $\mu$.

To this end, let $\mu^\dagger$ be the map
$$ \mu^\dagger: E \to \Hom(E, E) $$
adjoint to $\mu$. Consider the following diagram:

\begin{equation}\label{diag:leibniz_adjoint}
    \begin{tikzcd}
    	{E^A} && {\Hom(E^B,E^{A+B}) \oplus \Hom(E^{B'}, E^{(A+B)'})} \\
    	\\
    	{E^{A'}} && {\Hom(E^B,E^{(A+B)'})}
    	\arrow["{d^A}", from=1-1, to=3-1]
    	\arrow["{\mu^\dagger}", from=3-1, to=3-3]
    	\arrow["{\begin{bmatrix} \mu^\dagger \amsamp \mu^\dagger \end{bmatrix}}", from=1-1, to=1-3]
    	\arrow["{\begin{bmatrix} (d^{A+B})_* \\ (d^B)^* \end{bmatrix}}", from=1-3, to=3-3]
    \end{tikzcd}
\end{equation}

By adjointness, Diagrams \eqref{diag:leibniz_normal} and \eqref{diag:leibniz_adjoint} are equivalent.

We have an analogue of Proposition~\ref{prop:leibniz_easy}.
\begin{prop} \label{prop:leibniz_hard}
Let $D^B$ and $D^{A+B}$ be affine subspaces of $\Hom(E^B, E^{B'})$ and $\Hom(E^{A+B}, E^{(A+B)'})$ respectively. If the differential $d^B$ belongs to the subspace $D^B$ and the differential $d^{A+B}$ belongs to the subspace $D^{A+B}$, then the differential $d^A$ belongs to the affine subspace
$$ T(A,B) = (\mu^\dagger_*)^{-1}((\mu^\dagger)^*(D^{A+B}_*) + (\mu^\dagger)^*((D^B)^*)) $$
of $\Hom(E^{A}, E^{A'})$. Here we have
\begin{align*}
    D^{A+B}_* &\subseteq \Hom(\Hom(E^B, E^{A+B}), \Hom(E^B, E^{A+B'})), \\
    (\mu^\dagger)^*(D^{A+B}_*) &\subseteq \Hom(E^A, \Hom(E^B, E^{(A+B)'})), \\
    (D^B)^* &\subseteq \Hom(\Hom(E^{B'}, E^{A+B'}), \Hom(E^B, E^{A+B'})), \\
    (\mu^\dagger)^*((D^B)^*) &\subseteq \Hom(E^A, \Hom(E^B, E^{(A+B)'})).
\end{align*}
\end{prop}
\begin{proof}
As in the proof of Proposition~\ref{prop:leibniz_easy}, vector addition is a linear map. Furthermore, recall that the two maps
$$ (-)_*: \Hom(V, W) \to \Hom(\Hom(U, V), \Hom(U, W)), $$
$$ (-)^*: \Hom(V, W) \to \Hom(\Hom(W, U), \Hom(V, U)) $$
are also linear. With this notation, we can write
$$ T(A,B) = (\mu^\dagger_*)^{-1}(+((\mu^\dagger)^*((-)_*(D^{A+B})) \oplus (\mu^\dagger)^*((-)^*(D^B)))). $$
Repeated application of Lemma~\ref{lem:basic_affine} to this description shows that $T(A,B)$ is indeed an affine subspace of $\Hom(E^A, E^{A'})$. 

To see that it contains $d^A$, we note that equating the two composites in \eqref{diag:leibniz_adjoint} gives
$$ \mu^\dagger \circ d^{A+B}_* + \mu^\dagger \circ (d^B)^* = d^A \circ \mu^\dagger. $$
By assumption, this map is contained in the space
$$ (\mu^\dagger)^*(D^{A+B}_*) + (\mu^\dagger)^*((D^B)^*). $$
Taking preimages under $\mu^\dagger_*$ gives the desired result.
\end{proof}
\begin{rmk}
Proposition~\ref{prop:leibniz_hard} also holds when interchanging $A$ and $B$ throughout.
\end{rmk}

\section{The Algorithm}\label{sec:algo}

The pseudocode for our algorithm is laid out in Algorithm~\ref{alg:main}. A step-by-step explanation follows.

\begin{figure}
    \centering
    \begin{algorithm}[H]
        \SetAlgoLined
        \DontPrintSemicolon
        \For{$A$ a bidegree}{
            $D^A \longleftarrow \Hom(E^A, E^{A'})$\;
        }
        changed $\longleftarrow$ \texttt{true}\;
        \While{changed}{
            changed $\longleftarrow$ \texttt{false}\;
            \For{$A$ a bidegree}{
                \For{$B$ a bidegree, $B \geq A$}{
                    Compute $T(A,B)$\;
                    \If{$\dim(D^A \cap T(A,B)) < \dim(D^A)$}{
                        $D^A \longleftarrow D^A \cap T(A,B)$\;
                        changed $\longleftarrow$ \texttt{true}\;
                    }
                    Compute $T(B,A)$\;
                    \If{$\dim(D^B \cap T(B,A)) < \dim(D^B)$}{
                        $D^B \longleftarrow D^B \cap T(B,A)$\;
                        changed $\longleftarrow$ \texttt{true}\;
                    }
                    Compute $S(A,B)$\;
                    \If{$\dim(D^{A+B} \cap S(A,B)) < \dim(D^{A+B})$}{
                        $D^{A+B} \longleftarrow D^{A+B} \cap S(A,B)$\;
                        changed $\longleftarrow$ \texttt{true}\;
                    }
                }
            }
        }
        \caption{Leibniz propagation}\label{alg:main}
    \end{algorithm}
\end{figure}

First we initialize each space $D^A$ of differentials to be the whole space $\Hom(E^A, E^{A'})$. Recall that $D^A$ is an \emph{affine} subspace of $\Hom(E^A, E^{A'})$. We encode this data as a pair $v^A + U^A$ where $U^A$ is a linear subspace of the ambient space, and $v^A$ is a vector.

Next, we initialize a boolean flag \texttt{changed} to \texttt{true}. This flag will keep track of whether our algorithm has successfully reduced the dimension of some space $D^A$. Since we iterate over pairs of bidegrees in lexicographic order, it might happen that processing a pair of bidegrees reduces the dimension of some space $D^A$ which has already been processed. Example~\ref{ex:necessary_outer_loop} shows why this outer loop is necessary.

Then, we iterate over all pairs of bidegrees $A \leq B$ in the differential bigraded algebra $E$. The next three steps form the core of the algorithm.

By Propositions~\ref{prop:leibniz_easy} and \ref{prop:leibniz_hard}, the differentials on $E^A$, $E^B$, and $E^{A+B}$ must be contained in $T(A,B)$, $T(B,A)$, and $S(A,B)$ respectively. Therefore, we compute these spaces and intersect them with our current values of $D^A$, $D^B$ and $D^{A+B}$ respectively to get a potentially smaller subspace. Since the true differentials belong to both affine subspaces, their intersection is nonempty, so it is an affine subspace by Lemma~\ref{lem:affine_inter}.

\begin{example} \label{ex:necessary_outer_loop}
In the first iteration of the outer loop, when $A = (58,6)$ and $B = (76,16)$,
our algorithm makes the deduction that the diferential on the element $D_2$ must be either $h_0 Q_2$ or 0. Since this is a nontrivial deduction, it makes sure that the outer loop is run at least once more.

During the second iteration, when $A = (0,1)$ and $B = (58,6)$,
our algorithm uses our knowledge of the possible differentials on $D_2$ to infer that the differential on $h_0 D_2$ must be either $h_0^2 Q_2$ or 0. It was necessary to execute the outer loop a second time, because the pair $((0,1), (58,6))$ comes before $((58,6), (76,16))$ in the lexicographic order, so it had already been processed the first time (with no effect).
\end{example}

The following argument for the computation of $d_2(h_1) = 0$ is very well-known. It is also the first result given by our algorithm, which was found in only a few milliseconds.

\begin{example}
Consider $h_0$ and $h_1$ in bidegrees $(0,1)$ and $(1,1)$ respectively. We have that $h_0 h_1 = 0$, so in particular $d_2(h_0 h_1) = 0$. Moreover, $d_2(h_0) = 0$ for degree reasons, since the $E_2$ page has no non-zero elements in bidegree $(-1,3)$. Therefore,
$$ 0 = d_2(h_0 h_1) = d_2(h_0) h_1 + h_0 d_2(h_1) = h_0 d_2(h_1). $$
Since $d_2(h_1)$ has bidegree $(0,3)$, and since the only non-zero element $h_0^3$ in this bidegree is not $h_0$-torsion, we can conclude that $d_2(h_1) = 0$.
\end{example}

\bibliography{main}
\bibliographystyle{alphaurl}

\end{document}